\documentclass{amsart}

\usepackage{amssymb}
\usepackage[all]{xy}
\usepackage{hyperref}

\usepackage{enumitem}   

\setlist[enumerate]{itemsep=.2em,topsep=.2em,leftmargin=1.25em,itemindent=2.0em}

\newtheorem{thm}{Theorem}

\newtheorem{prop}[thm]{Proposition}

\theoremstyle{definition}
\newtheorem{defn}[thm]{Definition}

\newtheorem{say}[thm]{}
\newtheorem{exmp}[thm]{Example}

\newtheorem{rem}[thm]{Remark}          

\newtheorem*{ack}{Acknowledgments}      

\newtheorem{defn-thm}[thm]{Definition--Theorem}  
\newtheorem{defn-lem}[thm]{Definition--Lemma}  

\theoremstyle{remark}

\renewcommand{\c}[0]{{\mathbb C}}  

\renewcommand{\o}[0]{{\mathcal O}} 

\renewcommand{\r}[0]{{\mathbb R}}

\newcommand{\hh}[0]{{\mathbb H}}
\newcommand{\dd}[0]{{\mathbb D}}

\newcommand{\p}[0]{{\mathbb P}}

\newcommand{\qtq}[1]{\quad\mbox{#1}\quad}

\newcommand{\tsum}[0]{\textstyle{\sum}}

\def\loccoh#1.#2.#3.#4.{H^{#1}_{#2}(#3,#4)}

\DeclareMathAlphabet{\mathchanc}{OT1}{pzc}%
                                {m}{it}

\usepackage[all]{xy}\xyoption{dvips}

\newcommand{\tprod}[0]{\textstyle{\prod}}

\begin{document}
\bibliographystyle{amsalpha}

\title[Bounded meromorphic functions]{Bounded meromorphic functions \\ on the complex 2-disc}
  \author{J\'anos Koll\'ar}

\begin{abstract} 
  We describe  bounded, holomorphic functions   on the complex 2-disc, that admit   meromorphic extension  to a larger 2-disc.  This solves a conjecture of
  Bickel, Knese, Pascoe and Sola.
  The key technical ingredient is an old theorem of Zariski  about integrally closed ideals in 2-dimensional regular rings.
\end{abstract}

\maketitle

The aim of this note is to study bounded, holomorphic functions  $\psi$ on the 2-disc  $\dd^2:=\{(x,y): |x|,|y|< 1\}\subset \c^2$, that admit a 
meromorphic extension $\Psi$ to a larger 2-disc  $\{(x,y): |x|,|y|< 1+\epsilon\}$.

It is easy to see that the polar set of $\Psi$ intersects the closed 2-disc  $\bar \dd^2$ in finitely many points only, all of which satisfy $|x|=|y|=1$.
In Theorem~\ref{bkps.c.1.3.thm} we give a  complete local description of $\Psi$ at  these special points, answering  \cite[Conj.1.3]{bkps}.

We refer to \cite{bps, bkps} for a history of this question and to many related results. The higher dimensional cases are  of considerable interest, but appear  more complicated, see Remark~\ref{3d.rem}.

Up to conformal equivalence, this problem is the same as studying functions $\Phi(x, y)$ that are meromorphic in a neighborhood of the origin $(0,0)\in \c^2$ and bounded on
$\hh^2:=\{(x,y): \Im(x),\Im(y)>0\}\subset \c^2$.
From now on we write  $\hh^2_{\rm loc}$ 
to denote the intersection of $\hh^2$ with a suitable 
neighborhood of the origin $(0,0)\in \c^2$, and similarly for $\r^2_{\rm loc}, \c^2_{\rm loc}$.
We can write
$\Phi=f/g$ where $f, g\in \c\{x,y\}$ are analytic on $\c^2_{\rm loc}$ and have no common factors.

The first necessary condition is that $g(x,y)$  vanishes on $\overline{\hh}^2_{\rm loc}$ (the closure of $\hh^2_{\rm loc}$)  only at the origin.
A complete description of such analytic functions is given in \cite{bkps}, in terms of their Puiseux factors as follows.

By Newton, we can write any $g\in \c\{x,y\}$ (that is not divisible by $x$) uniquely as a product of  Puiseux series
$$
g=u(x,y)\tprod_i \bigl(y+\phi_i(x^{1/r_i})\bigr),
$$
where $u(0,0)\neq 0$, the $\phi_i$ are holomorphic, $\phi_i(0)=0$ and $r_i\geq 1$.

\begin{thm}\cite[p.4]{bkps}\label{bkps.p.fact.thm} The holomorphic function $g(x,y)\in \c\{x,y\}$  vanishes on $\overline{\hh}^2_{\rm loc}$ only at the origin  iff each of its  Puiseux  factors  is of the form
  $$
  y+q_i(x)+ x^{m_i}\psi_i(x^{1/r_i}),
  \eqno{(\ref{bkps.p.fact.thm}.1)}
  $$
  where   $q_i\in\r\{x\}$, $q_i(0)=0$, $q'_i(0)>0$,  the $m_i$ are even and $\Im (\psi_i(0))>0$.
\end{thm}
This gives a complete description of the possible denominators; see Proposition~\ref{how.to.check} on how to check this condition in concrete situations.

In order to understand the numerators,   \cite[p.7]{bkps}
defines ideals as follows.

\begin{defn}\label{ideals.I.defn} For the Puiseux  factors in 
Theorem~\ref{bkps.p.fact.thm} set
$$
{\mathcal I}\bigl(y+q_i(x)+ x^{m_i}\psi_i(x^{1/r_i})\bigr):=\bigl(y+q_i(x), x^{m_i}\bigr).
$$
Then define ${\mathcal I}(g)$ as   their product
$$
 {\mathcal I}(g):=\tprod_i \bigl(y+q_i(x), x^{m_i}\bigr)\subset  \c\{x,y\}.
  $$
\end{defn}

The following  theorem
 gives a positive answer to \cite[Conj.1.3]{bkps}.

\begin{thm} \label{bkps.c.1.3.thm}
  Let $\Phi=f(x,y)/g(x,y)$ be a meromorphic function on $\c^2_{\rm loc}$,  where $f, g\in \c\{x,y\}$ are analytic on $\c^2_{\rm loc}$ and have no common factors.
  Then $\Phi$ is bounded on  $\hh^2_{\rm loc}$ iff
  \begin{enumerate}
    \item 
  $g(x,y)$  vanishes on $\overline{\hh}^2_{\rm loc}$ only at the origin, and
    \item $f(x, y)\in {\mathcal I}(g)$.
      \end{enumerate}
\end{thm}

 The conditions are sufficient by \cite[Thm.1.2]{bkps}, which also shows the necessity if $g$ has either multiplicity 2 or an ordinary singularity.  The necessity of (\ref{bkps.c.1.3.thm}.1) 
 follows from  Theorem~\ref{bkps.p.fact.thm}, the new claim is that (\ref{bkps.c.1.3.thm}.2) is also necessary.

\begin{rem}  
 The answer of Theorem~\ref{bkps.c.1.3.thm} seems asymmetrical in $x,y$, but this is an accident of our choices.

  If $G(x,y)=ax+by+\cdots$ with $ab\neq 0$ then, by the Weierstrass preparation theorem, for every $m>1$ there are unique polynomials $q_1(x), q_2(y)$ of degree $<m$ such that  $q_i(0)=0$, and
$$
\bigl(G, (x,y)^m\bigr)=\bigl(y+q_1(x), x^m\bigr)=\bigl(x+q_2(y), y^m\bigr).
$$
\end{rem}

\begin{exmp} \cite[2.18]{bkps}  The simplest non-trivial examples come from 
  $$
  (y+x+ix^{2m})^2-x^{4m+1} =(y+x+ix^{2m}+x^{2m+1/2})(y+x+ix^{2m}-x^{2m+1/2}).
  $$
  The corresponding ideal is $(y+x,x^{2m})^2=\bigl((y+x)^2, (y+x)x^{2m}, x^{4m}\bigr)$.
  Thus we get that meromorphic functions with denominator $ (y+x+ix^{2m})^2-x^{4m+1}$, that are bounded on $\hh^2_{\rm loc}$ are of the form
  $$
  \frac{v_1(x,y)(y+x)^2+ v_2(x,y) (y+x)x^{2m} +v_3(x,y) x^{4m}}{(y+x+ix^{2m})^2-x^{4m+1}},
  $$
  where the $v_i(x,y)\in\c\{x,y\}$ are holomorphic.
\end{exmp}

In the first step of the proof, we transform the boundedness of $f/g$  on $\hh^2_{\rm loc}$ into an inequality over $\c^2_{\rm loc}$.

\begin{thm}\label{4.EQUIV.THM}
  Using the notation of Theorems~\ref{bkps.p.fact.thm}~and~\ref{bkps.c.1.3.thm}, the following are equivalent.
  \begin{enumerate}
  \item $f/g$ is bounded on $\hh^2_{\rm loc}$.
  \item $|f|\leq   C\cdot \tprod_i\bigl(|y+q_i(x)|+|x^{m_i}|\bigr)$ for $(x, y)\in \hh^2_{\rm loc}$.
    \item $|f|\leq   C\cdot \tprod_i\bigl(|y+q_i(x)|+|x^{m_i}|\bigr)$ for $(x, y)\in \r^2_{\rm loc}$.
    \item $|f|\leq   C\cdot \tprod_i\bigl(|y+q_i(x)|+|x^{m_i}|\bigr)$ for $(x, y)\in \c^2_{\rm loc}$.
  \end{enumerate}
\end{thm}

Then we use the theory of integral closure of ideals and a theorem of  \cite{zar38} to show that (\ref{4.EQUIV.THM}.4) is equivalent to (\ref{bkps.c.1.3.thm}.2).

\begin{say}[Proof of Theorem~\ref{4.EQUIV.THM}; beginning]\label{pf.step.1.say}
  The implications  (\ref{4.EQUIV.THM}.2) $\Rightarrow$ (\ref{4.EQUIV.THM}.3) and (\ref{4.EQUIV.THM}.4) $\Rightarrow$ (\ref{4.EQUIV.THM}.2) are clear.

  The equivalence (\ref{4.EQUIV.THM}.1) $\Leftrightarrow$ (\ref{4.EQUIV.THM}.2)
   is  just a reformulation of
  \cite[Thm.1.2]{bkps}, which shows that, for every $i$,
  $$
|y+q_i(x)|+|x^{m_i}|\leq   C\cdot |y+\phi_i(x)| \qtq{for} (x, y)\in \hh^2_{\rm loc}.
$$
Since the reverse bound is clear, we get that
$|g(x, y)|$ and $\tprod_i\bigl(|y+q_i(x)|+|x^{m_i}|\bigr)$ are mutually bounded by constant multiples of each other. Thus the conditions
$$
|f(x,y)|\leq  C\cdot |g(x,y)|\qtq{and} |f(x,y)|\leq  C\cdot \tprod_i\bigl(|y+q_i(x)|+|x^{m_i}|\bigr)
$$
are equivalent on $\hh^2_{\rm loc}$. \qed
\end{say}

The advantage of the form (\ref{4.EQUIV.THM}.4) is that it ties in with the notion of integral dependence of holomorphic functions.

\begin{defn}[Integral dependence]   Let $X$ be a complex space, $x\in X$ a point and
$\o_{X, x}$ the ring of germs of holomorphic functions at $x$.
  Let $J=(g_1,\dots, g_m)\subset \o_{X, x}$ be an ideal  and
  $h\in  \o_{X, x}$. Then $h$ is {\it integral} over $J$ if $|h|\leq C\cdot\sum_i |g_i|$ for some $C>0$ in some neighborhood of $x\in X$.
  This notion is   independent of the generators chosen.

  $J$ is called {\it integrally closed} if it contains every holomorphic function that is integral over $J$.

  Most of the literature on integral dependence is  algebraic, with a very different definition. Among the standard books, the 
 equivalence is  stated in  \cite[Rem.9.6.10]{laz-book} and \cite[Thm.7.1.7]{swa-hun},  with further references for proofs.
\end{defn}

We need the valuative criterion of integral dependence.
However, we need not only its statement as in \cite{laz-book, swa-hun}, but a more precise version  treated in \cite{MR2499856}.

\begin{say}[Valuative criterion of integral dependence]\cite{MR2499856} \label{val.crit.dep.say}
  Let $J=(g_1, \dots, g_r)\subset \o_{X,x}$ be an ideal sheaf and $f$ a holomorphic fucntion. How do we check that $f$ is integral over $J$?

  First, if $X=\c^n$ and $J=(\prod_i x_i^{m_i})$ is  generated by a monomial, then clearly
  $f$ is integral over $J$ iff it vanishes along $(x_i=0)$ with multiplicity
  $\geq m_i$ for every $i$. 
  This can be checked along $n$ arcs as follows.
  
  For each $i$ let  $\gamma_i:\dd\to \c^n$ be an analytic arc
  whose image intersects the hyperplane $(x_i=0)$ at a single ploint $\gamma_i(0)$, and is disjoint from the $(x_j=0)$ for $j\neq i$.
  Then  $f$ is integral over $J=(g:=\prod_i x_i^{m_i})$ iff the $(f/g)\circ \gamma_i$ are all bounded near $0\in \dd$.
  Note also that since $f/g $ is meromorphic, it is enough to check boundedness on a set that accumulates to $0\in \dd$.

  Most ideals are not monomial, but Hironaka's resolution theorem says
  that for any reduced complex space $X$ and ideal sheaf $J\subset \o_X$,
  there is a proper, bimeromorphic morphism  $\pi:Y\to X$
  such that $Y$ is smooth and the preimage  $J_Y:=\pi^*J\subset \o_{Y}$ is locally monomial. This is called a {\it monomialization} of $J$.

  Combining with our previous discussions, we get the {\it valuative criterion:}
  \medskip

  {\it Theorem \ref{val.crit.dep.say}.1.} $f$ is integral over $J=(g_1, \dots, g_r)\subset \o_X$ iff, for every analytic arc $\gamma:\dd\to X$, there is a $C_{\gamma}>0$ such that
  $$
  |f\circ\gamma(t)|\leq C_{\gamma}\cdot \tsum_i |g_i\circ\gamma(t)|.
  \eqno{(\ref{val.crit.dep.say}.1.1)}
  \qed
  $$
  \medskip

  Note, however, that our previous discussions give a  more precise result.
  Let $E_j\subset Y$ be the $\pi$-exceptional divisors.
Choose arcs    $\gamma_j:\dd\to Y$ such that 
$\gamma_j(\dd)$ intersects  $E_j$ at a single point $\gamma_j(0)$ and is disjoint from the other $E_k$ for $k\neq j$.

\medskip

{\it Refined Theorem \ref{val.crit.dep.say}.2.} $f$ is integral over $J=(g_1, \dots, g_r)\subset \o_{X,x}$ iff, for the above analytic arcs $\pi\circ\gamma_j:\dd\to X$,
there are $C_j>0$ and sets $U_j\subset \dd$ accumulating to $0\in \dd$, such that 
  $$
  |f\circ\gamma_j(t)|\leq C_j\cdot \tsum_i |g_i\circ\gamma_j(t)| \qtq{for $t\in U_j$.}
  \eqno{(\ref{val.crit.dep.say}.2.1)}
  $$
  
  \medskip
\end{say}

Note that the ideals  ${\mathcal I}(g)$ in Definition~\ref{ideals.I.defn} are generated by real power series. This turns out to be quite important, so we need some general facts about real power series and real structures on complex manifolds.

\begin{say}[Reality questions]\label{real.p.s.say}
  Let $h\in \r\{x,y\}$  be an irreducible power series. It either stays 
  irreducible in $\c\{x,y\}$, or it decomposes as a product of 2 irreducible factors. A typical example is $y^2+x^2$.

  Let $C_h:=(h=0)\subset \c^2_{\rm loc}$ be the corresponding curve.
  Then $h$ is irreducible in $\c\{x,y\}$ iff the normalization
  $\bar C_h\to C_h$ has a unique point $c\in \bar C_h$ lying over the origin.
  Complex conjugation of $\c^2$ lifts to an
  antiholomorphic involution in $ \bar C_h$  and $c$ is a fixed point.

  If $p\in \c^2$ is a real point, then complex conjugation  lifts to the blow-up  $B_p\c^2$, giving a real structure.
The exceptional curve of $B_p\c^2\to \c^2$ is then isomorphic to $\c\p^1$ with its standard real structure  $(x{:}y)\mapsto (\bar x{:}\bar y)$.
The same applies after any number of blow-ups of real points  $\pi:Y\to \c^2$. 
\end{say}

\begin{prop}\label{pf.step.2.prop}
 Let  $J=(g_1,\dots, g_m)\subset \r\{x,y\}$ be a real ideal and  $f\in  \c\{x,y\}$ an analytic function.  Assume that $J$  has a monomialization $\pi:Y\to\c^2$ where we blow up only real points.
  Then the following are equivalent.
  \begin{enumerate}
  \item $|f(x,y)|\leq C\cdot \sum_i |g_i(x,y)|$ for all $(x,y)\in \r^2_{\rm loc}$ for some $C>0$.
     \item $|f(x,y)|\leq C\cdot \sum_i |g_i(x,y)|$ for all $(x,y)\in \c^2_{\rm loc}$ for some $C>0$.
  \end{enumerate}
\end{prop}

Proof. The implication (\ref{pf.step.2.prop}.2) $\Rightarrow$ (\ref{pf.step.2.prop}.1) is clear.

To see the converse,
we follow the approach explained in (\ref{val.crit.dep.say}.2).
By assumption $Y$ has a real structure compatible with the standard one on $\c^2$, and   all the $\pi$-exceptional curves $E_j\subset Y$ are  real isomorphic to $\c\p^1$.
Thus, for every $E_j$  we can choose a real analytic arc $\gamma_j:\dd\to Y$ that intersects $E_j$ transversally  but is disjoint from the other exceptional curves.
Thus $\pi\circ\gamma_j$ maps $\dd\cap \r$ to $\r^2$.

If (\ref{pf.step.2.prop}.1) holds then
(\ref{val.crit.dep.say}.2.1) is satisfied with $U_j:=\dd\cap \r$.
Thus $f$ is integral over $J$ by Theorem~\ref{val.crit.dep.say}.2,
which is (\ref{pf.step.2.prop}.2),
 \qed

\begin{say}[Monomialization  for ideals in $\c\{x,y\}$]\label{monids.say}  
While the general monomialization method of Hironaka is quite complicated,  there is a  simple algorithm for ideals $J=(g_1,\dots, g_m)\subset \c\{x,y\}$.

Start with the curve  $C=(g_1\cdots g_m=0)$ and $Y_0:=\c^2$.
If $\pi_i:Y_i\to \c^2$ is already constructed, then we get $Y_{i+1}\to Y_i$ by blowing up all the  intersection points  of the
$\pi_i$-exceptional curve and of  the birational transform $C_i\subset Y_i$ of $C$; that is, the closure of
$\pi_i^{-1}\bigl(C\setminus\{(0,0)\}\bigr)$.
After finitely many steps, the preimage of the ideal $J$ becomes monomial.
(This algorithm depends on the choice of the $g_i$, but works well for our purposes. See \cite[Chap.1]{k-res} for a discussion of various other methods.)

Assume now that the $g_i$ are real and
write  $\prod_ig_i=\prod_j p_{j}$ as a product of irreducible factors  in $\r\{x,y\}$.
Then  we repeatedly blow up  the
intersection points  of the exceptional curve and the union of  the birational transforms 
 of the curves $B_{j}:=(p_{j}=0)$.  If the $p_{j}$ are also
irreducible   in $\c\{x,y\}$, then, as we noted in Paragraph~\ref{real.p.s.say}, for each $B_j$ the intersection consists of a single real point. Then
at each step we blow up only real points. We have thus proved the following.
\medskip

{\it Claim \ref{monids.say}.1.}  Let  $J=(g_1,\dots, g_m)\subset \r\{x,y\}$
be an ideal. Assume that every  $\r\{x,y\}$-irreducible factor of $\prod_i g_i$ 
is irreducible in $\c\{x,y\}$. Then  $J$  has a monomialization $\pi:Y\to\c^2$ where we blow up only real points.\qed
\medskip

\end{say}

\begin{say}[Proof of Theorem~\ref{4.EQUIV.THM}; end]\label{pf.step.2.say}
  We prove the last implication  (\ref{4.EQUIV.THM}.3) $\Rightarrow$ (\ref{4.EQUIV.THM}.4).
     
  Our ideal 
  ${\mathcal I}(g)$  is the product of ideals of the form
  $\bigl(y+q(x), x^m\bigr)$.
  Thus it has a generating set   $g_j$ where each $g_j$ is a product of a power of $x$ with some
  $\prod_{\ell} (y+q_{\ell}(x))$. Thus the irreducible factors of $\prod_jg_j$ are $x$ and the $y+q_i(x)$. These are irreducible both in 
  $\r\{x,y\}$ and $\c\{x,y\}$.

  Thus ${\mathcal I}(g)$  has a monomialization $\pi:Y\to\c^2$ where we blow up only real points by Claim~\ref{monids.say}.1.
   Therefore Proposition~\ref{pf.step.2.prop}  applies to  ${\mathcal I}(g)$,  and
  shows that  (\ref{4.EQUIV.THM}.3) implies (\ref{4.EQUIV.THM}.4). \qed
   \end{say}

\begin{say}[Proof of Theorem~\ref{bkps.c.1.3.thm}]\label{pf.step.3.say}
  If we expand $\tprod_i\bigl(|y+q_i(x)|+|x^{m_i}|\bigr)$, we get the sum of the absolute values of the generators of the ideal  ${\mathcal I}(g)$.
  Thus  Theorem~\ref{4.EQUIV.THM}.4 says that 
 $f$ is  integral over
  ${\mathcal I}(g)$. 
  If ${\mathcal I}(g)=\tprod_i\bigl(y+q_i(x), x^{m_i}\bigr)$ is integrally closed, then   $f\in {\mathcal I}(g)$, which is what we want.

  Note that  $\bigl(y+q_i(x), x^{m_i}\bigr)$  is, after a coordinate change, the same as  $(y, x^{m_i})$, which is clearly integrally closed.
  By a theorem of Zariski, in a two-dimensional regular local ring, products of integrally closed ideals are integrally closed.
  This is proved in \cite{zar38}; a  more accessible reference is
  \cite[Thm.14.4.4]{swa-hun}.

  Since $\c\{x,y\}$ is two-dimensional and regular,
  we obtain that  $\tprod_i\bigl(y+q_i(x), x^{m_i}\bigr)$ is integrally closed. Therefore 
  $f\in \tprod_i\bigl(y+q_i(x), x^{m_i}\bigr)={\mathcal I}(g)$. \qed
\end{say}

\begin{rem}\label{3d.rem}
  Zariski's theorem  holds only in dimension 2. The following example, based on \cite{bps} and e-mails of A.~Sola and I.~Swanson, shows that  higher dimensional versions of Theorem~\ref{bkps.c.1.3.thm} are more complicated.

The ideals
$(z+x+y, x^2+y^2)$ and $(z+x+y, x^2+2y^2)$ are integrally closed, but their
product
$$
\bigl((z+x+y)^2, (z+x+y)x^2, (z+x+y)y^2, (x^2+y^2)(x^2+2y^2)\bigr)
$$
is not since $(z+x+y)xy$ is integral over it.
This gives that
$$
\frac{(z+x+y)xy}
     {(z+x+y+ix^2+iy^2)(z+x+y+ix^2+2iy^2)}
     $$
     is   bounded on $\hh^3_{loc}$, but it is not a linear combination of products of bounded functions with denominators
     $(z+x+y+ix^2+iy^2)$ and $(z+x+y+ix^2+2iy^2)$.
\end{rem}
     


  The conclusion of Theorem~\ref{bkps.p.fact.thm} can be checked without writing down the Puiseux series, working with 1 irreducible factor $p(x,y)$ of $g$ at a time.
 Newton's method of rotating rulers (see, for example,  \cite[pp.126-127]{newt-II} or \cite[Thm.18]{arcology} for a detailed treatment) constructs the Puiseux series solution of $p(x,y)=0$ term by term as
  $$
  y=a_1x+a_2x^2+\cdots+a_mx^{m}+\cdots.
  $$
  The formulas become more complicated once a fractional power of $x$ appears, but the method is very transparent while we have only integer powers.
We obtain the following.

  \begin{prop} \label{how.to.check}
  Let $p(x,y)\in \c\{x,y\}$. Assume that $p(0,0)=0$ and the lowest $y$ power in $p$ is $y^r$ (with coefficient 1). Let $h(x)\in \c[x]$ be a polynomial of degree $m$ such that $h(0)=0$. Then the following hold.
  \begin{enumerate}
  \item If $p$ is irreducible and one Puiseux factor of $p$  has the form
    $$
    y+h(x)+(\mbox{higher terms}),
    \eqno{(\ref{how.to.check}.1.1)}
    $$ then so is every Puiseux factor.
    \item Every Puiseux factor of $p$ has the form
      (\ref{how.to.check}.1.1) iff we can write
      $$
      p(x, y)=(y+h(x))^r+\tsum_{i+j/m>r}(y+h(x))^ix^j.
      $$
    \item The implicit function theorem solves  
$$
      \frac1{r!}\frac{\partial^{r-1} p}{\partial y^{r-1}}=0\qtq{as}
      y=-h(x)-(\mbox{higher terms}).
      $$
    \item If $p$ is a Weierstrass polynomial in $y$ (of degree $r$) then
      $$
      \frac1{r!}\frac{\partial^{r-1} p}{\partial y^{r-1}}=y+h(x)+(\mbox{higher terms}).
      $$
      \end{enumerate}
  \end{prop}

  Note, however, that in (\ref{how.to.check}.1.1) the higher terms form a Puiseux  series (thus involve fractional powers of $x$), but in (\ref{how.to.check}.3--4) the higher terms form a power  series.

  When we apply this to Theorem~\ref{bkps.p.fact.thm}, we use
  $h(x)=q(x)+\psi(0)x^m$.

  \begin{ack} I thank
K.~Bickel, G.~Knese, J.~E.~Pascoe and A.~Sola for the Conjecture and  e-mails;
    L.~Lempert and I.~Swanson for helpful suggestions and  references.
Partial  financial support    was provided  by  the NSF under grant number
DMS-1901855.
\end{ack}


\def\cprime{$'$} \def\cprime{$'$} \def\cprime{$'$} \def\cprime{$'$}
  \def\cprime{$'$} \def\dbar{\leavevmode\hbox to 0pt{\hskip.2ex
  \accent"16\hss}d} \def\cprime{$'$} \def\cprime{$'$}
  \def\polhk#1{\setbox0=\hbox{#1}{\ooalign{\hidewidth
  \lower1.5ex\hbox{`}\hidewidth\crcr\unhbox0}}} \def\cprime{$'$}
  \def\cprime{$'$} \def\cprime{$'$} \def\cprime{$'$}
  \def\polhk#1{\setbox0=\hbox{#1}{\ooalign{\hidewidth
  \lower1.5ex\hbox{`}\hidewidth\crcr\unhbox0}}} \def\cdprime{$''$}
  \def\cprime{$'$} \def\cprime{$'$} \def\cprime{$'$} \def\cprime{$'$}
\providecommand{\bysame}{\leavevmode\hbox to3em{\hrulefill}\thinspace}
\providecommand{\MR}{\relax\ifhmode\unskip\space\fi MR }
\providecommand{\MRhref}[2]{%
  \href{http://www.ams.org/mathscinet-getitem?mr=#1}{#2}
}
\providecommand{\href}[2]{#2}

\bigskip

  Princeton University, Princeton NJ 08544-1000, \

\email{kollar@math.princeton.edu}

\end{document}